\pdfoutput=1
\RequirePackage{ifpdf}
\ifpdf 
\documentclass[pdftex]{sigma}
\else
\documentclass{sigma}
\fi

\usepackage{enumitem}

\numberwithin{equation}{section}

\newtheorem{Theorem}{Theorem}[section]
\newtheorem*{Theorem*}{Theorem}
\newtheorem{Corollary}[Theorem]{Corollary}
\newtheorem{Lemma}[Theorem]{Lemma}
\newtheorem{Proposition}[Theorem]{Proposition}

\theoremstyle{definition}
\newtheorem{Definition}[Theorem]{Definition}

\newtheorem{Remark}[Theorem]{Remark}

\DeclareMathOperator{\im}{im}

\DeclareMathOperator{\sfl}{sf}
\DeclareMathOperator{\diag}{diag}

\DeclareMathOperator{\GL}{GL}

\DeclareMathOperator{\codim}{codim}
\DeclareMathOperator{\gra}{graph}

\begin{document}

\allowdisplaybreaks

\newcommand{\arXivNumber}{2508.01061}

\renewcommand{\thefootnote}{}

\renewcommand{\PaperNumber}{029}

\FirstPageHeading

\ShortArticleName{On the Uniqueness of the $G$-Equivariant Spectral Flow}

\ArticleName{On the Uniqueness of the $\boldsymbol{G}$-Equivariant Spectral Flow\footnote{This paper is a~contribution to the Special Issue on Asymptotics, Randomness and Noncommutativity. The~full collection is available at \href{https://sigma-journal.com/noncommutativity.html}{https://sigma-journal.com/noncommutativity.html}}}

\Author{Marek IZYDOREK~$^{\rm a}$, Joanna JANCZEWSKA~$^{\rm a}$, Maciej STAROSTKA~$^{\rm b}$\newline and Nils WATERSTRAAT~$^{\rm b}$}

\AuthorNameForHeading{M.~Izydorek, J.~Janczewska, M.~Starostka and N.~Waterstraat}

\Address{$^{\rm a)}$~Institute of Applied Mathematics, Faculty of Applied Physics and Mathematics, \\
\hphantom{$^{\rm a)}$}~Gda\'{n}sk University of Technology,
Narutowicza 11/12, 80-233 Gda\'{n}sk, Poland}
\EmailD{\mail{marek.izydorek@pg.edu.pl}, \mail{joanna.janczewska@pg.edu.pl}}

\Address{$^{\rm b)}$~Martin-Luther-Universit\"at Halle-Wittenberg, Naturwissenschaftliche Fakult\"at II,\\
\hphantom{$^{\rm b)}$}~Institut f\"ur Mathematik, 06099 Halle (Saale), Germany}
\EmailD{\mail{maciej.starostka@pg.edu.pl}, \mail{nils.waterstraat@mathematik.uni-halle.de}}

\ArticleDates{Received August 06, 2025, in final form March 02, 2026; Published online March 25, 2026}

\Abstract{The spectral flow is an integer-valued homotopy invariant for paths of selfadjoint Fredholm operators. Lesch as well as Pejsachowicz, Fitzpatrick and Ciriza independently showed that it is uniquely characterised by its elementary properties. The authors recently introduced a $G$-equivariant spectral flow for paths of selfadjoint Fredholm operators that are equivariant under the action of a compact Lie group $G$. The purpose of this paper is to show that the $G$-equivariant spectral flow is uniquely characterised by the same elementary properties when appropriately restated. As an application, we introduce an alternative definition of the $G$-equivariant spectral flow via a $G$-equivariant Maslov index.}


\Classification{58J30; 47A53; 37C81}

\renewcommand{\thefootnote}{\arabic{footnote}}
\setcounter{footnote}{0}

\section{Introduction}

Let $H$ be a real separable Hilbert space of infinite dimension and let us denote by $\mathcal{L}(H)$ the bounded, by $\GL(H)$ the invertible, by $\mathcal{K}(H)$ the compact and by $\mathcal{S}(H)$ the selfadjoint operators on $H$. An operator $T\in\mathcal{L}(H)$ is called Fredholm if its kernel and cokernel are of finite dimension. For $T\in\mathcal{S}(H)\subset\mathcal{L}(H)$, these conditions are equivalent to a finite-dimensional kernel and a closed range. Atiyah and Singer considered in~\cite{AtiyahSinger} the set $\mathcal{FS}(H)\subset\mathcal{L}(H)$ of all selfadjoint Fredholm operators and showed that it consists of the three connected components
\begin{align*}
\mathcal{FS}^+(H)&:=\{T\in\mathcal{FS}(H)\mid \sigma_{\rm ess}(T)\subset(0,+\infty)\},\\
\mathcal{FS}^-(H)&:=\{T\in\mathcal{FS}(H)\mid \sigma_{\rm ess}(T)\subset(-\infty,0)\},
\end{align*}
and
\begin{align*}
\mathcal{FS}^i(H):=\mathcal{FS}(H)\setminus\mathcal{FS}^\pm(H).
\end{align*}
Here $\sigma_{\rm ess}(T)$ denotes the essential spectrum of $T$, i.e., the set of all $\lambda\in\mathbb{R}$ such that $\lambda-T$ is not a Fredholm operator. Operators in $\mathcal{FS}^\pm(H)$ have as only possible points in the spectrum below/above $0$ eigenvalues of finite multiplicity. Instead operators in $\mathcal{FS}^i(H)$ necessarily have essential spectrum above and below~$0$.

Atiyah and Singer also showed in \cite{AtiyahSinger} that the two components $\mathcal{FS}^\pm(H)$ of $\mathcal{FS}(H)$ are contractible as topological spaces, whereas $\mathcal{FS}^i(H)$ has a rather rich topology and in particular an infinitely cyclic fundamental group. In a joint work with Patodi \cite{AtiyahPatodi}, they introduced an explicit isomorphism between $\pi_1(\mathcal{FS}^i(H))$ and the integers. Actually, the idea of this map, which is called spectral flow, is sensible for any closed or open path and in any of the three components of $\mathcal{FS}(H)$. Following Phillips' analytic approach \cite{Phillips} (cf.~\cite{book}), its construction can be outlined as follows: Let $I=[0,1]$ be the unit interval and $L=\{L_\lambda\}_{\lambda\in I}$ a path in $\mathcal{FS}(H)$. As the spectra of the operators $L_\lambda$ cannot accumulate at $0$ (cf. \cite{ProcHan,Fredholm}), there is a partition $0=\lambda_0<\dots< \lambda_N=1$ of the unit interval and numbers $a_i>0$, $i=1,\dots, N$, such that the maps
\begin{align}\label{specproj}
[\lambda_{i-1},\lambda_i]\ni \lambda\mapsto\chi_{[-a_i,a_i]}(L_\lambda)\in\mathcal{L}(H)
\end{align}
are continuous families of finite rank projections. Here $\chi_{[a,b]}(T)$ denotes the spectral projection of a selfadjoint operator $T$ with respect to the interval $[a,b]\subset\mathbb{R}$. Consequently, the spaces $E(L_\lambda,[0,a_i]):=\im(\chi_{[0,a_i]}(L_\lambda))$, $\lambda_{i-1}\leq \lambda\leq \lambda_i$, $i=1,\dots,N$, are of finite dimension and the spectral flow of the path $L$ can be defined by
\begin{align}\label{sfl}
\sfl(L)=\sum^N_{i=1}{(\dim(E(L_{\lambda_i},[0,a_i]))-\dim(E(L_{\lambda_{i-1}},[0,a_i])))}\in\mathbb{Z}.
\end{align}
Note that $\dim(E(L_\lambda,[0,a]))$ is the number of eigenvalues of $L_\lambda$ in the interval $[0,a]$ including multiplicities. Roughly speaking, the spectral flow is the net number of eigenvalues of $L_0$ crossing~$0$ from the negative to the positive half of the axis whilst the parameter is travelling along the interval. It can be shown that \eqref{sfl} is well defined, i.e., it only depends on the path~$L$ and neither on the choice of the partition $\{\lambda_i\}_{i=0,\dots,N}$ of $I$ nor on the numbers $a_1,\dots,a_N$ as long as the maps~\eqref{specproj} are continuous and finite-rank valued. The spectral flow has the following properties:
\begin{itemize}\itemsep=0pt
 \item[$(\mathcal{Z})$] If $L=\{L_\lambda\}_{\lambda\in I}$ is a path in $\mathcal{FS}(H)\cap\GL(H)$, then
\[
 \sfl(L)=0.
\]
 \item[$(\mathcal{C})$] If $L^1$, $L^2$ are paths in $\mathcal{FS}(H)$, then
 \[\sfl\bigl(L^1\ast L^2\bigr)=\sfl\bigl(L^1\bigr)+\sfl\bigl(L^2\bigr),\]
 where $L^1\ast L^2$ denotes the concatenation of $L^1$ and $L^2$.
 \item[$(\mathcal{A})$] If $H=H_1\oplus H_2$ and $L_\lambda(H_i)\subset H_i$, $i=1,2$, for a path $L=\{L_\lambda\}_{\lambda\in I}$ in $\mathcal{FS}(H)$, then
 \[\sfl(L)=\sfl(L|_{H_1})+\sfl(L|_{H_2}).\]
 \item[$(\mathcal{H})$] If $h \colon I\times I\rightarrow\mathcal{FS}(H)$ is a homotopy such that $h(s,0), h(s,1)\in\GL(H)$ for all $s\in I$, then
 \[\sfl(h(0,\cdot))=\sfl(h(1,\cdot)).\]
\end{itemize}
Let us note that $(\mathcal{Z})$, $(\mathcal{C})$, $(\mathcal{A})$ easily follow from the definition~\eqref{sfl}, and $(\mathcal{H})$ is a simple modification of the original argument in~\cite{Phillips} (see, e.g.,~\cite{CompSfl}).

It is a natural question to ask whether some of these properties uniquely characterise the spectral flow when a normalisation property is added that in particular excludes the trivial map which assigns $0\in\mathbb{Z}$ to any path in $\mathcal{FS}(H)$. Henceforth, we denote for a pair of topological spaces $(X,Y)$ by $\Omega(X,Y)$ the set of all paths in $X$ having endpoints in $Y$. Lesch considered in \cite{Lesch} maps $\mu \colon \Omega(\mathcal{FS}^i(H),\mathcal{FS}^i(H)\cap\GL(H))\rightarrow\mathbb{Z}$ and proved that any such map which satisfies the same properties $(\mathcal{C})$, $(\mathcal{H})$ as the spectral flow is already given by \eqref{sfl} if the following normalisation assumption holds:
\begin{itemize}\itemsep=0pt
\item[$\mathcal{(N)}$] There is some $T_0\in\mathcal{FS}^i(H)\cap\GL(H)$ and some rank-one orthogonal projection $P$ in $H$ commuting with $T_0$ such that
\[ \mu\bigl(\{\lambda P+(I_H-P)T_0\}_{-\frac{1}{2}\leq\lambda\leq\frac{1}{2}} \bigr)=1 \]
\end{itemize}
Note that the path \smash{$\{\lambda P+(I_H-P)T_0\}{}_{\smash{-\frac{1}{2}\leq\lambda\leq\frac{1}{2}}}$} is in $\GL(H)$ for any $\lambda\neq 0$, whereas $0$ is a simple eigenvalue for $\lambda=0$. It is readily seen from \eqref{sfl} that indeed \smash{$\sfl(\{\lambda P\hspace{-0.2pt}+(I_H\hspace{-0.2pt}-P)T_0\}{}_{\smash{-\frac{1}{2}\leq\lambda\leq\frac{1}{2}}} )=1$}.

Another approach to the uniqueness problem of the spectral flow was proposed by Fitzpatrick, Pejsachowicz and Ciriza in \cite{JacoboUniqueness}. They considered maps $\mu \colon \Omega(\mathcal{FS}(H),\mathcal{FS}(H)\cap\GL(H))\rightarrow\mathbb{Z}$ and showed that any such map which satisfies the same properties $(\mathcal{Z})$, $(\mathcal{A})$ and $(\mathcal{H})$ as the spectral flow is already given by \eqref{sfl} if the following normalisation assumption holds:
\begin{itemize}\itemsep=0pt
\item[$\mathcal{(M)}$] If $H$ is of finite dimension and $L=\{L_\lambda\}_{\lambda\in I}$ is a path in $(\mathcal{FS}(H),\mathcal{FS}(H)\cap\GL(H))$, then
\[\mu(L)=\mu_{\rm Morse}(L_0)-\mu_{\rm Morse}(L_1),\]
where $\mu_{\rm Morse}(L_\lambda)$ denotes the Morse index of $L_\lambda$, $\lambda\in I$, i.e., the number of negative eigenvalues including multiplicities.
\end{itemize}
Inspired by~\cite{Fang}, the authors defined in \cite{MJN21} a $G$-equivariant spectral flow for paths $L=\{L_\lambda\}_{\lambda\in I}$ in $\mathcal{FS}(H)$ where each operator $L_\lambda$ is equivariant under the orthogonal action of a compact Lie group $G$ on $H$. The $G$-equivariant spectral flow is an element of the real representation ring $\mathrm{RO}(G)$ and it formally has the same properties $(\mathcal{Z})$, $(\mathcal{C})$, $(\mathcal{A})$, $(\mathcal{H})$ and $(\mathcal{M})$ when they are appropriately restated in terms of virtual representations of $G$. The aim of this paper is to show that the $G$-equivariant spectral flow is uniquely characterised by these properties as well. If the group action is trivial, then our theorem is the uniqueness of the spectral flow \eqref{sfl} from \cite{JacoboUniqueness}.

We recall the definition of the $G$-equivariant spectral flow in the next section, where we also introduce the announced uniqueness of it as our main theorem. Section 3 is devoted to the rather elaborate proof of the latter theorem. In the final Section 4, we introduce a $G$-equivariant Maslov index, where we adopt an approach to the classical Maslov index from a celebrated paper by Cappell, Lee and Miller \cite{Cappell}. Finally, we show as first application of our uniqueness theorem that the $G$-equivariant spectral flow of a path $L=\{L_\lambda\}_{\lambda\in I}$ can be obtained as $G$-equivariant Maslov index of the graphs of the operators $L_\lambda$, which is well-known in the non-equivariant case~\cite{CompSfl}.

\section{The \texorpdfstring{$\boldsymbol{G}$}{G}-equivariant spectral flow and its uniqueness}\label{sect-gequivsfl}

Let $(V,\rho)$ be a real representation of the compact Lie group $G$, i.e., $V$ is a finite-dimensional real vector space and $\rho \colon G\rightarrow\GL(V)$ is a group homomorphism. Two representations $(V_1,\rho_1)$, $(V_2,\rho_2)$ of $G$ are isomorphic if there is a $G$-equivariant isomorphism $\alpha \colon V_1\rightarrow V_2$, i.e., $\rho_2(g)\circ \alpha=\alpha\circ\rho_1(g)$ for all $g\in G$. The set of isomorphism classes of representations of $G$ is a commutative monoid when representations are added by the direct sum. Elements of the associated Grothendieck group $\mathrm{RO}(G)$ are formal differences $[U]-[V]$ of isomorphism classes of $G$-representations modulo the equivalence relation generated by $[U]-[V]\sim [U\oplus W]-[V\oplus W]$. Consequently, the inverse element of $[U]-[V]$ is $[V]-[U]$ and the neutral element in $\mathrm{RO}(G)$ is $[V]-[V]$ for any $G$-representation~$V$ (see~\cite{Segal}).

Let now $G$ be a compact Lie group that acts orthogonally on the real separable Hilbert space $H$. In what follows, we denote by $\mathcal{FS}(H)^G$ the set of $G$-equivariant selfadjoint Fredholm operators, i.e.,
\[T(gu)=g(T u),\qquad u\in H,\quad g\in G,\]
and we let $\mathcal{FS}^\pm(H)^G$ and $\mathcal{FS}^i(H)^G$ be the corresponding subsets of the connected components of $\mathcal{FS}(H)$.

Let now $L=\{L_\lambda\}_{\lambda\in I}$ be a path in $\mathcal{FS}(H)^G$, i.e., a continuous map from the unit interval $I$ into the space $\mathcal{FS}(H)^G$. Since the operators $L_\lambda$ are $G$-equivariant, it follows that the spaces $E(L_{\lambda},[0,a])$ in \eqref{sfl} are $G$-invariant. Consequently they define equivalence classes of $G$-representations and so the idea of \eqref{sfl} carries over to $\mathrm{RO}(G)$ by setting
\begin{align}\label{sfl-equiv}
\sfl_G(L)=\sum^N_{i=1}{([E(L_{\lambda_i},[0,a_i])]-[E(L_{\lambda_{i-1}},[0,a_i])])}\in \mathrm{RO}(G).
\end{align}
Following Phillips argument for \eqref{sfl} from \cite{Phillips}, the authors showed in \cite{MJN21} that \eqref{sfl-equiv} only depends on the path $L$ and neither on the choice of the partition $0=\lambda_0<\dots< \lambda_N=1$ of the unit interval nor on the numbers $a_i>0$, $i=1,\dots,N$. Thus, this equivariant spectral flow is well defined. Note that if $G=\{e\}$ is the trivial group, then representations are isomorphic if and only if they are of the same dimension. Thus, $\mathrm{RO}(\{e\})\cong\mathbb{Z}$ and \eqref{sfl-equiv} can be identified with the ordinary spectral flow~\eqref{sfl}. Moreover, there is a canonical homomorphism
\begin{align*}
F\colon \ \mathrm{RO}(G)\rightarrow\mathbb{Z},\qquad [U]-[V]\mapsto\dim(U)-\dim(V),
\end{align*}
and we see from \eqref{sfl} and \eqref{sfl-equiv} that
\begin{align*}
F(\sfl_G(L))=\sfl(L).
\end{align*}
Consequently, if $\sfl_G(L)$ is trivial, then the classical spectral flow of $L$ has to vanish as well. On the other hand, it is not difficult to find a simple example of a path of $G=\mathbb{Z}_2$-equivariant operators such that $\sfl(L)=0\in\mathbb{Z}$ but $\sfl_G(L)\in \mathrm{RO}(\mathbb{Z}_2)\cong\mathbb{Z}\oplus\mathbb{Z}$ is non-trivial (cf.\ \cite[Section~2.4]{MJN21}). Examples for continuous groups can be found in~\cite{MJN25}. Other references for $G$-equivariant spectral flows are, e.g., \cite{Fang, Hochs,Lim,Liu}.

We now consider, for a given compact Lie group, maps
\[\mu=\mu_{H}\colon \ \Omega\bigl(\mathcal{FS}(H)^G,\GL(H)\cap\mathcal{FS}(H)^G\bigr)\rightarrow \mathrm{RO}(G)\]
that are supposed to exist for any separable Hilbert space $H$. Such maps might have the following properties, where as before $L=\{L_\lambda\}_{\lambda\in I}$ is a path of bounded operators on $H$:
\begin{itemize}[leftmargin=30pt]\itemsep=0pt
 \item[$(\mathcal{Z}_G)$] If $L_\lambda\in\GL(H)\cap\mathcal{FS}(H)^G$ for all $\lambda\in I$, then
 \begin{align*}
 \mu(L)=0\in \mathrm{RO}(G).
 \end{align*}
 \item[$(\mathcal{A}_G)$] Let $H=H_1\oplus H_2$, where $H_1$, $H_2$ are $G$-invariant, and let $L=\{L_\lambda\}_{\lambda\in I}$ be a path with invertible endpoints in $\mathcal{FS}(H)^G$ such that $L_\lambda(H_i)\subset H_i$ for $i=1,2$, and all $\lambda\in I$. Then
\begin{align*}
\mu(L)=\mu(L|_{H_1})+\mu(L|_{H_2})\in \mathrm{RO}(G).
\end{align*}
\item[$(\mathcal{H}_G)$] If $h \colon I\times I\rightarrow\mathcal{FS}(H)^G$ is a homotopy such that $h(s,0)$ and $h(s,1)$ are invertible for all $s\in I$, then
\begin{align*}
\mu(h(0,\cdot))=\mu(h(1,\cdot))\in \mathrm{RO}(G).
\end{align*}
\item[$(\mathcal{M}_G)$] If $H$ is of finite dimension and $L=\{L_\lambda\}_{\lambda\in I}$ is a path in $\bigl(\mathcal{FS}(H)^G,\mathcal{FS}(H)^G\cap\GL(H)\bigr)$, then
\[\mu(L)=[E^-(L_0)]-[E^-(L_1)]\in \mathrm{RO}(G),\]
where $E^-(T)$ denotes the direct sum of the eigenspaces with respect to negative eigenvalues of a selfadjoint operator $T \colon H\rightarrow H$.
\end{itemize}
The $G$-equivariant spectral flow \eqref{sfl-equiv} as a map
\[\sfl_G \colon \ \Omega\bigl(\mathcal{FS}(H)^G,\GL(H)\cap\mathcal{FS}(H)^G\bigr)\rightarrow \mathrm{RO}(G)\]
has all these properties. This was shown for $(\mathcal{Z}_G)$ in \cite[Lemma~2.5]{MJN21}, for $(\mathcal{A}_G)$ in \cite[Lemma~2.6]{MJN21}, for $(\mathcal{H}_G)$ in \cite[Corollary~2.10]{MJN21} and for $(\mathcal{M}_G)$ in \cite[Proposition~3.2]{MJN21}. Now we can state the main theorem of this paper.

\begin{Theorem}\label{thm-uniqueness}
Let $G$ be a compact Lie group and
\[\mu=\mu_{H} \colon \ \Omega\bigl(\mathcal{FS}(H)^G,\GL(H)\cap\mathcal{FS}(H)^G\bigr)\rightarrow \mathrm{RO}(G)\]
a map that is defined for any separable Hilbert space $H$. If $\mu$ satisfies $(\mathcal{Z}_G)$, $(\mathcal{A}_G)$, $(\mathcal{H}_G)$ and~$(\mathcal{M}_G)$, then
\[\mu=\sfl_G.\]
\end{Theorem}

If $G$ is the trivial group, we can identify \eqref{sfl} and \eqref{sfl-equiv} to obtain the main theorem of Fitzpatrick, Pejsachowicz and Ciriza's work~\cite{JacoboUniqueness}.

\begin{Corollary}\label{cor-uniqueness}
Let
\[\widetilde{\mu}=\widetilde{\mu}_{H}\colon \ \Omega(\mathcal{FS}(H),\GL(H)\cap\mathcal{FS}(H))\rightarrow\mathbb{Z}\]
be a map that is defined for any separable Hilbert space $H$. If $\widetilde{\mu}$ satisfies $(\mathcal{Z})$, $(\mathcal{A})$, $(\mathcal{H})$ and~$(\mathcal{M})$, then
\[\widetilde{\mu}=\sfl.\]
\end{Corollary}

Note that \eqref{sfl-equiv} can be verbatim defined on complex Hilbert spaces as an element of the complex representation ring~$R(G)$. The properties $(\mathcal{Z}_G)$, $(\mathcal{A}_G)$, $(\mathcal{H}_G)$, $(\mathcal{M}_G)$ and Theorem~\ref{thm-uniqueness} still hold in this case, but one step in the proof of Theorem~\ref{thm-uniqueness} becomes redundant, as we explain below in Section \ref{sect-OG}.

Finally, let us note that there also is a corresponding concatenation property:
\begin{itemize}\itemsep=0pt
 \item[$(\mathcal{C}_G)$] If $L^1,L^2\in\Omega\bigl(\mathcal{FS}(H)^G,\GL(H)\cap\mathcal{FS}(H)^G\bigr)$, then
 \[\mu(L^1\ast L^2)=\mu\bigl(L^1\bigr)+\mu\bigl(L^2\bigr)\in \mathrm{RO}(G),\]
\end{itemize}
which clearly holds for the $G$-equivariant spectral flow by its definition \eqref{sfl-equiv}. Actually, by arguing as in \cite[Proposition~4.26]{Robbin-Salamon}, any map $\mu\colon \Omega\bigl(\mathcal{FS}(H)^G,\GL(H)\cap\mathcal{FS}(H)^G\bigr)\rightarrow \mathrm{RO}(G)$ that satisfies $(\mathcal{Z}_G)$, $(\mathcal{A}_G)$ and $(\mathcal{H}_G)$ has the property $(\mathcal{C}_G)$. On the other hand, it is readily seen that $(\mathcal{H}_G)$ and $(\mathcal{C}_G)$ imply $(\mathcal{Z}_G)$ (cf.\ \cite[Lemma~5.3]{Lesch}) and thus the axioms are strongly related. Note that, however, $(\mathcal{C}_G)$ cannot replace $(\mathcal{A}_G)$ in Theorem~\ref{thm-uniqueness}, as the map $\mu$ which is defined by $(\mathcal{M}_G)$ if $\dim H<\infty$, and $0\in \mathrm{RO}(G)$ otherwise, would satisfy all axioms. Thus, $(\mathcal{C}_G)$~as characterising property of the $G$-equivariant spectral flow requires a different normalisation property than~$(\mathcal{M}_G)$. Recall that Lesch's uniqueness theorem~\cite{Lesch} for the classical spectral flow~\eqref{sfl} assumes $(\mathcal{N})$ instead of $(\mathcal{M})$. Of course, it is a natural question whether some appropriate normalisation condition $(\mathcal{N}_G)$ in combination with $(\mathcal{C}_G)$ and $(\mathcal{H}_G)$ also characterises the $G$-equivariant spectral flow. Let us note that Lesch's argument in \cite[Section~5.3]{Lesch} uses the connectedness of some subspaces of selfadjoint Fredholm operators which in general does not hold in the $G$-equivariant case (cf., e.g., Step~2 in the proof of Theorem~\ref{thm-cogredpar} below) and thus~\cite{Lesch} apparently does not transfer without major modifications.

\section{Proof of Theorem~\ref{thm-uniqueness}}

\subsection{The \texorpdfstring{$\boldsymbol{G}$}{G}-equivariant cogredient parametrix}
Henceforth, we denote by $\mathcal{KS}(H)$ the space of all selfadjoint compact operators with the norm topology, and by $\mathcal{KS}(H)^G$ its $G$-equivariant elements. Similarly, $\GL(H)^G$ are the $G$-equivariant invertible operators. The aim of this section is to introduce a transformation of paths in $\mathcal{FS}(H)^G$ into a normal form. To be more precise, we show that for any path $L=\{L_\lambda\}_{\lambda\in I}$ in $\mathcal{FS}(H)^G$, there are a path $\{M_\lambda\}_{\lambda\in I}$ in $\GL(H)^G$, a path $\{K_\lambda\}_{\lambda\in I}$ in $\mathcal{KS}(H)^G$ and a selfadjoint operator $Q\in\GL(H)^G$, such that
\[M^\ast_\lambda L_\lambda M_\lambda=Q+K_\lambda,\qquad\lambda\in I.\]
The path $\{M_\lambda\}_{\lambda\in I}$ is called $G$-equivariant cogredient parametrix for $L=\{L_\lambda\}_{\lambda\in I}$ and the particular form of the operator $Q$ will be specified below. Before we begin the construction, let us recall that if $T\in\mathcal{L}(H)$ is a $G$-equivariant selfadjoint operator and $f \colon \sigma(T)\rightarrow\mathbb{R}$ is continuous, then the selfadjoint operator $f(T)$ is $G$-equivariant as well (see, e.g., \cite[Lemma~3.3]{MJN21}). In~what follows, we will use this fact without further reference. We now firstly consider paths in $\mathcal{FS}^\pm(H)^G$.

\begin{Theorem}\label{thm-cogrpm}
Every path $L=\{L_\lambda\}_{\lambda\in I}$ in $\mathcal{FS}^\pm(H)^G$ has a cogredient parametrix. More precisely, there is a path $\{M_\lambda\}_{\lambda\in I}$ in $\GL(H)^G$ and a path $\{K_\lambda\}_{\lambda\in I}$ in $\mathcal{KS}(H)^G$ such that
\[M^\ast_\lambda L_\lambda M_\lambda=\pm I_H+K_\lambda,\qquad\lambda\in I.\]
\end{Theorem}

\begin{proof}
We first note that it suffices to show the claim for paths in $\mathcal{FS}^+(H)^G$ as otherwise we can consider the path $-L=\{-L_\lambda\}_{\lambda\in I}$.

The proof of the theorem is based on the following simple fact. If $T\in\mathcal{FS}^+(H)^G$ and $\chi_{(0,\infty)}(T)$, $\chi_{(-\infty,0]}(T)$ are the spectral projections, then $\chi_{(0,\infty)}(T)+\chi_{(-\infty,0]}(T)=I_H$ and we obtain
\begin{align*}
T&=\bigl(\chi_{(0,\infty)}(T)+\chi_{(-\infty,0]}(T)\bigr)T\bigl(\chi_{(0,\infty)}(T)+\chi_{(-\infty,0]}(T)\bigr)\\
&=\chi_{(0,\infty)}(T)T\chi_{(0,\infty)}(T)+\chi_{(-\infty,0]}(T)T\chi_{(-\infty,0]}(T)\\
&=\bigl(\chi_{(0,\infty)}(T)T\chi_{(0,\infty)}(T)+\chi_{(-\infty,0]}(T)\bigr)+\bigl(\chi_{(-\infty,0]}(T)T\chi_{(-\infty,0]}(T)-\chi_{(-\infty,0]}(T)\bigr)\\
&=:S+K,
\end{align*}
where $S\in\GL(H)^G$ is positive and $K\in\mathcal{KS}(H)^G$. Thus, for any $\lambda_0\in I$, we can decompose~$L_{\lambda_0}$ as $L_{\lambda_0}=S_{\lambda_0}+K_{\lambda_0}$, where $S_{\lambda_0}\in\GL(H)^G$ is positive and $K_{\lambda_0}\in\mathcal{KS}(H)^G$. As $\GL(H)$ is open, there is some interval $I'\subset I$ containing $\lambda_0$ such that $S'_\lambda:=L_\lambda-K_{\lambda_0}\in\GL(H)^G$ for any $\lambda\in I'$. By the continuity of spectra, $S'_\lambda$ actually is positive for all $\lambda\in I'$. Consequently, if we set $K'_\lambda:=K_{\lambda_0}$, we have $L_\lambda=S'_\lambda+K'_\lambda$ for $\lambda\in I'$, where $S'_\lambda\in\GL(H)^G$ is positive and $K'_\lambda\in\mathcal{KS}(H)^G$.

By compactness, there is a finite cover of $I$ by intervals $I_1,\dots, I_k$, and positive operators $S^i_\lambda\in\GL(H)^G$ as well as $K^i_\lambda\in\mathcal{KS}(H)^G$ such that $L_\lambda=S^i_\lambda+K^i_\lambda$ for $\lambda\in I_i$, $i=1,\dots,k$. We~take a partition of unity $\{\eta_i\}_{1\leq i\leq k}$ subordinate to the covering $\{I_i\}_{1\leq i\leq k}$ of $I$ and set for $\lambda\in I$
\[S_\lambda:=\sum^k_{i=1}{\eta_i(\lambda)\, S^i_\lambda},\qquad K_\lambda:=\sum^k_{i=1}{\eta_i(\lambda)\, K^i_\lambda}.\]
Since the set of positive $G$-equivariant operators as well as the set of $G$-equivariant selfadjoint compact operators are convex, taking sums and multiplying by the maps $\{\eta_i\}_{1\leq i\leq k}$ preserve these properties. Thus, we finally conclude that $L_\lambda=S_\lambda+K_\lambda$ for $\lambda\in I$, where $S_\lambda\in\GL(H)^G$ is positive and $K_\lambda\in\mathcal{KS}(H)^G$.

As $S_\lambda$ is positive, $M_\lambda:=\bigl(S^\frac{1}{2}_\lambda\bigr)^{-1}\in\GL(H)^G$ is indeed defined and we see that
\[M^\ast_\lambda L_\lambda M_\lambda=M^\ast_\lambda S_\lambda M_\lambda+M^\ast_\lambda K_\lambda M_\lambda=I_H+M^\ast_\lambda K_\lambda M_\lambda,\]
which shows the theorem since $M^\ast_\lambda K_\lambda M_\lambda\in\mathcal{KS}(H)^G$.
\end{proof}

For paths in $\mathcal{FS}^i(H)^G$, the argument is far more sophisticated and a detailed discussion can be found in the authors' recent work \cite{MJN25}. Here we include a concise version of the proof for the convenience of the reader. An operator $Q\in\mathcal{L}(H)$ is called a symmetry if it is selfadjoint and orthogonal, which in particular implies that $Q^2=I_H$. Any symmetry is of the form $Q=P-(I_H-P)=2P-I_H$ for an orthogonal projection $P$, and $Q\in\mathcal{FS}^i(H)$ if and only if $P$ has infinite-dimensional kernel and range. Finally, $Q$ is $G$-equivariant if and only if the subspaces $\im(P)$ and $\ker(P)$ of $H$ are $G$-invariant.

\begin{Theorem}\label{thm-cogredpar}
Every path $L=\{L_\lambda\}_{\lambda\in I}$ in $\mathcal{FS}^i(H)^G$ has a cogredient parametrix. More precisely, there is a $G$-equivariant symmetry $Q\in\mathcal{FS}^i(H)^G$ as well as paths $M=\{M_\lambda\}_{\lambda\in I}$ in $\GL(H)^G$ and $K=\{K_\lambda\}_{\lambda\in I}$ in $\mathcal{KS}(H)^G$ such that
\[M^\ast_\lambda L_\lambda M_\lambda=Q+K_\lambda,\qquad\lambda\in I.\]
\end{Theorem}

\begin{proof}
We split the proof into three steps.

{\it Step~$1$: The map $\pi_Q$ and local sections.}
We consider for $G$-equivariant symmetries $Q\in\mathcal{FS}^i(H)^G$ the maps
\begin{align}\label{defQ}
\pi_Q \colon \ \GL(H)^G\times\mathcal{KS}(H)^G\rightarrow\mathcal{FS}^i(H)^G,\qquad \pi_Q(M,K)=MQM^\ast+K.
\end{align}
Note that Theorem~\ref{thm-cogredpar} is proved if we can show that there is some $Q$ for which there is a~continuous map $\tilde{L} \colon I\rightarrow \GL(H)^G\times\mathcal{KS}(H)^G$ such that $L_\lambda=\pi_Q\circ\tilde{L}$ for all $\lambda\in I$, i.e., the path $L$ can be lifted to $\GL(H)^G\times\mathcal{KS}(H)^G$ with respect to $\pi_Q$. The right setting to find such a lifting is the theory of fibre bundles and our first aim is to construct local sections for \eqref{defQ}. To be more precise, we are going to show that for any $S\in\mathcal{FS}^i(H)^G$ there is a $G$-equivariant symmetry $Q_S$, an open neighbourhood $\mathcal{U}_S$ of $S$ in $\mathcal{FS}^i(H)^G$ and a map $\sigma_S \colon \mathcal{U}_S\rightarrow\GL(H)^G\times\mathcal{KS}(H)^G$ such that
\[(\pi_{Q_S}\circ\sigma_S)(T)=T\qquad\text{for all} \ T\in\mathcal{U}_S.\]
To construct $\sigma_S$, let $K$ be the orthogonal projection onto the kernel of the selfadjoint Fredholm operator $S$. Then $\ker(S)$ is $G$-invariant and of finite dimension, and thus $K\in\mathcal{KS}(H)^G$ as well~as
\begin{align}\label{V}
V:=S+K\in \GL(H)^G.
\end{align}
To simplify notation, we henceforth let $P_+(V)=\chi_{(0,\infty)}(V)$ and $P_-(V)=\chi_{(-\infty,0)}(V)$ be the orthogonal projections onto the positive and negative spectral subspaces of $V$, which are $G$\nobreakdash-equiv\-ari\-ant. The operator
\begin{align}\label{QS}
Q_S=2P_+(V)-I_H\in\mathcal{FS}^i(H)^G
\end{align}
has a neighbourhood $\tilde{\mathcal{U}}\subset\mathcal{FS}^i(H)^G$ of invertible operators. As before, we denote for $T\in\tilde{\mathcal{U}}$ by $P_+(T)$ and $P_-(T)$ the orthogonal projections onto the positive and negative spectral subspaces. As they continuously depend on $T\in\tilde{\mathcal{U}}$, $Q_S$ has a neighbourhood $\mathcal{U}\subset\tilde{\mathcal{U}}$ such that
\[P_+(T) P_+(Q_S)+P_-(T) P_-(Q_S)\in\GL(H)^G,\qquad T\in \mathcal{U},\]
where it is used that this operator is $I_H\in\GL(H)^G$ for $T=Q_S$. This shows that the restriction of $P_\pm(T)$ to $\im(P_\pm(Q_S))$ is a bijection onto $\im(P_\pm(T))$ for any $T\in\mathcal{U}$. To simplify notation, we~set $H_\pm:=\im(P_\pm(Q_S))$ and consider for $T\in\mathcal{U}$ the map $B(T) \colon H_+\times H_+\rightarrow\mathbb{R}$,
\[B(T)(u,v)=\langle TP_+(T)u,P_+(T)v\rangle,\qquad u,v\in H_+.\]
Then $B(T)$ is positive definite, as well as $G$-invariant as $T$, $P_+(T)$ are $G$-equivariant and $G$ acts orthogonally. Let $T_B$ be the unique selfadjoint operator on $H_+$ such that
\[B(T)(u,v)=\langle T_Bu,v\rangle,\qquad u,v\in H_+,\]
which is usually called the Riesz-representation of $B(T)$. The selfadjoint operator $T_B$, and thus \smash{$S_+(T):=T_B^{\smash{-\frac{1}{2}}}$}, is $G$\nobreakdash-equiv\-ari\-ant as it is unique and $B(T)$ is $G$-invariant. Moreover,
\begin{align*}
\langle u,v\rangle&=\bigl\langle T_B^{-1}T_Bu,v\bigr\rangle=\bigl\langle T_B^{-\frac{1}{2}}T_Bu,T_B^{-\frac{1}{2}}v\bigr\rangle=\bigl\langle T_BT_B^{-\frac{1}{2}}u,T_B^{-\frac{1}{2}}v\bigr\rangle=B(T)\bigl(T_B^{-\frac{1}{2}}u,T_B^{-\frac{1}{2}}v\bigr)\\
&=\bigl\langle TP_+(T)T_B^{-\frac{1}{2}}u,P_+(T)T_B^{-\frac{1}{2}}v\bigr\rangle=\bigl\langle T_B^{-\frac{1}{2}}P_+(T)TP_+(T)T_B^{-\frac{1}{2}}u,v\bigr\rangle\\
&=\langle S_+(T)P_+(T)TP_+(T)S_+(T)u,v\rangle,\qquad u,v\in H_+,
\end{align*}
which implies that
\begin{align}\label{equ-S+}
S_+(T)P_+(T)TP_+(T)S_+(T)=I_{H_+}.
\end{align}
Analogously, there is a map $S_- \colon \mathcal{U}\rightarrow\GL(H_-)^G$ such that for all $T\in\mathcal{U}$
\begin{align}\label{equ-S-}
S_-(T)P_-(T) TP_-(T)S_-(T)=-I_{H_-}.
\end{align}
We now define $S_0 \colon \mathcal{U}\rightarrow\GL(H)^G$ by
\[S_0(T)=P_+(T)S_+(T)P_+(Q_S)-P_-(T)S_-(T)P_-(Q_S),\]
and note that by \eqref{equ-S+} and \eqref{equ-S-}
\begin{align*}
S_0(T)^\ast TS_0(T)=P_+(Q_S)-P_-(Q_S)=Q_S,\qquad T\in\mathcal{U}.
\end{align*}
This shows that
\[\sigma \colon \ \mathcal{U}\rightarrow\GL(H)^G\times\mathcal{KS}(H)^G,\qquad \sigma(T)=((S_0(T)^{-1})^\ast,0),\]
satisfies
\begin{align}\label{equ-picircsigma}
(\pi_{Q_S}\circ\sigma)(T)=\pi_{Q_S}\bigl(\bigl(S_0(T)^{-1}\bigr)^\ast,0\bigr)=\bigl(S_0(T)^{-1}\bigr)^\ast Q_S S_0(T)^{-1}=T,\qquad T\in\mathcal{U}.
\end{align}
Let us recall that $\mathcal{U}$ is a neighbourhood of the symmetry $Q_S$ in \eqref{QS}. We now aim to construct a section $\sigma_V$ in a neighbourhood $\mathcal{U}_V$ of the map $V$ that was defined in \eqref{V}. The product $\GL(H)^G\times\mathcal{KS}(H)^G$ is a closed subgroup of the topological group $\GL(H)\times\mathcal{KS}(H)$, where
\begin{align}\label{groupmult}
(M,K)\cdot \bigl(\tilde{M},\tilde{K}\bigr)=\bigl(M\tilde{M}, K+M\tilde{K}M^\ast\bigr),\qquad (M,K), \bigl(\tilde{M},\tilde{K}\bigr)\in\GL(H)\times\mathcal{KS}(H).
\end{align}
Moreover, there is an action $\tau$ of $\GL(H)^G\times\mathcal{KS}(H)^G$ on $\mathcal{FS}^i(H)^G$ defined by
\[\tau_{h}(L)= MLM^\ast+K,\qquad L\in\mathcal{FS}^i(H),\qquad h=(M,K)\in\GL(H)^G\times\mathcal{KS}(H)^G.\]
For $h:=\bigl(|V|^\frac{1}{2},0\bigr)\in \GL(H)^G\times\mathcal{KS}(H)^G$, we obtain from \eqref{QS}
\[\tau_h(Q_S)=|V|^\frac{1}{2}Q_S|V|^\frac{1}{2}=V\]
and thus see that $\mathcal{U}_V:=\tau_h(\mathcal{U})$ is an open neighbourhood of $V$ in $\mathcal{FS}^i(H)^G$. We set
\[\sigma_V \colon \ \mathcal{U}_V\rightarrow \GL(H)^G\times\mathcal{KS}(H)^G,\qquad \sigma_V(T)=h\cdot \sigma(\tau_{h^{-1}}(T)),\qquad T\in\mathcal{U}_V,\]
where here and henceforth the dot stands for the group multiplication \eqref{groupmult}. To show that $\sigma_V$ is a section of $\pi_{Q_S}$ over $\mathcal{U}_V$, we first note that
\[\pi_{Q_S}(h_1\cdot h_2)=\tau_{h_1}(\pi_{Q_S}(h_2)),\qquad h_1,h_2\in \GL(H)^G\times\mathcal{KS}(H)^G\]
which is a consequence of the definition of $\tau$ and \eqref{groupmult}. Thus, it follows from \eqref{equ-picircsigma} that
\begin{align}
(\pi_{Q_S}\circ\sigma_V)(T)&=\pi_{Q_S}(h\cdot \sigma(\tau_{h^{-1}}(T)))=\tau_h(\pi_{Q_S}(\sigma(\tau_{h^{-1}}(T))))\nonumber\\
&=\tau_h(\tau_{h^{-1}}(T))=T,\qquad T\in\mathcal{U}_V,\label{QSsigmaV}
\end{align}
as claimed. Finally, we set $\tilde{h}=(I_H,-K)\in \GL(H)^G\times\mathcal{KS}(H)^G$, consider the open neighbourhood $\mathcal{U}_S=\tau_{\tilde{h}}(\mathcal{U}_V)$ of $S$ and define
\[\sigma_S\colon \ \mathcal{U}_S\rightarrow \GL(H)^G\times\mathcal{KS}(H)^G,\quad \sigma_S(T)=\tilde{h}\cdot \sigma_V(\tau_{\tilde{h}^{-1}}(T)),\qquad T\in\mathcal{U}_S.\]
The same computation as in \eqref{QSsigmaV} shows that $\pi_{Q_S}\circ\sigma_S=\mathrm{id}|_{\mathcal{U}_S}$, which finishes our first step.

{\it Step $2$: The image of $\pi_Q$.}
The aim of this step of the proof of Theorem~\ref{thm-cogredpar} is to show that for any symmetry $Q$, the image $\pi_Q$ is a union of connected components of $\mathcal{FS}^i(H)^G$. In the argument, we need the auxiliary result that for any symmetries $Q_1$ and $Q_2$
\begin{align}\label{lem33}
\im(\pi_{Q_1})\cap \im(\pi_{Q_2})\neq\varnothing \ \Longrightarrow \ \im(\pi_{Q_1})=\im(\pi_{Q_2})\subset\mathcal{FS}^i(H)^G.
\end{align}
We postpone the proof of \eqref{lem33} and first obtain the actual aim of this step by showing that $\im(\pi_Q)$ is open and closed in $\mathcal{FS}^i(H)^G$. By the first step of the proof, for any $S\in\im(\pi_Q)$ there is a symmetry $Q_S$ and an open neighbourhood $\mathcal{U}_S$ of $S$ in $\mathcal{FS}^i(H)^G$ such that $\mathcal{U}_S\subset\im(\pi_{Q_S})$. Now \eqref{lem33} implies that $\mathcal{U}_S\subset\im(\pi_Q)$ and thus $\im(\pi_Q)$ is open in $\mathcal{FS}^i(H)^G$. Secondly, let $\{S_n\}_{n\in\mathbb{N}}\subset\im(\pi_Q)$ converge to some $S\in\mathcal{FS}^i(H)^G$. Again, by the first step of the proof, there is a symmetry $Q_S$ and a neighbourhood $\mathcal{U}_S$ of $S$ in $\mathcal{FS}^i(H)^G$ such that $\mathcal{U}_S\subset\im(\pi_{Q_S})$. As~$\{S_n\}_{n\in\mathbb{N}}$ converges to $S$, we see that $S_n\in\mathcal{U}_S$ for sufficiently large $n$. By~\eqref{lem33}, this implies that $\mathcal{U}_S\subset\im(\pi_Q)$ and thus $S\in\im(\pi_Q)$, which means that $\im(\pi_Q)$ is closed in $\mathcal{FS}^i(H)^G$.

It remains to show \eqref{lem33}, for which we let $Q_1$, $Q_2$ be two symmetries and assume that $S\in \im(\pi_{Q_1})\cap \im(\pi_{Q_2})$. Then there are $(M_1,K_1), (M_2,K_2)\in\GL(H)^G\times\mathcal{KS}(H)^G$ such that
\[S=\pi_{Q_1}(M_1,K_1)=\pi_{Q_2}(M_2,K_2),\]
and thus
\[S=M_1Q_1M^\ast_1+K_1=M_2Q_2M^\ast_2+K_2,\]
which implies
\[Q_2=M^{-1}_2M_1Q_1M^\ast_1\bigl(M^{-1}_2\bigr)^\ast+M^{-1}_2(K_1-K_2)\bigl(M^{-1}_2\bigr)^\ast.\]
A direct computation shows that
\[\pi_{Q_2}(h)=\pi_{Q_1}\bigl(h\cdot \tilde{h}\bigr),\qquad h\in \GL(H)^G\times\mathcal{KS}(H)^G,\]
where $\tilde{h}=\bigl(M^{-1}_2M_1,M^{-1}_2(K_1-K_2)\bigl(M^{-1}_2\bigr)^\ast\bigr)\in \GL(H)^G\times\mathcal{KS}(H)^G$. This yields $\im(\pi_{Q_2})\subset\im(\pi_{Q_1})$, and proves \eqref{lem33} by swapping $Q_1$ and $Q_2$.

{\it Step $3$: Fibre bundles and end of the proof.}
Henceforth, we set $B_Q=\im(\pi_Q)$ for any given symmetry $Q\in\mathcal{FS}^i(H)^G$. Moreover, let us note that if $S=Q$ in Step 1, then it follows from \eqref{V} and \eqref{QS} that $Q_Q=Q$, and thus for every symmetry $Q\in\mathcal{FS}^i(H)^G$ there is some $S\in\mathcal{FS}^i(H)^G$ such that $Q_S=Q$. Therefore, it is not necessary to specify $S$ in our notation any longer.

We now first note that the map $\pi_Q \colon \GL(H)^G\times\mathcal{KS}(H)^G\rightarrow B_Q$ is the projection of a locally trivial fibre-bundle with fibre given by the isotropy group of $Q\in\mathcal{FS}^i(H)^G$. Indeed, by the first step of the proof, there is a local section $\sigma$ of $\pi_Q$ on an open subset $\mathcal{U}_Q\subset B_Q$. Now a local trivialisation over $\mathcal{U}_Q$ is given by
\[\eta \colon \ \mathcal{U}_Q\times\pi^{-1}_Q(Q)\rightarrow\pi^{-1}_Q(\mathcal{U}_Q),\qquad \eta(S,h)=\sigma(S)\cdot h,\]
and it can be shifted to any point $T\in B_Q$ in the following usual way. As $\im(\pi_Q)=B_Q$ by definition, there is some element $\tilde{h}\in \GL(H)^G\times\mathcal{KS}(H)^G$ such that $\pi_Q(\tilde{h})=\tau_{\tilde{h}}(Q)=T$. Now $\mathcal{U}:=\tau_{\tilde{h}}(\mathcal{U}_Q)$ is a neighbourhood of $T$ and $\tau_{\tilde{h}} \colon \mathcal{U}_Q\rightarrow\mathcal{U}$ is a homeomorphism, which yields a local trivialisation over $\mathcal{U}$ by
\[\eta' \colon \ \mathcal{U}\times\pi^{-1}_Q(Q)\rightarrow\pi^{-1}_Q(\mathcal{U}),\qquad \eta'(S,h)=\tilde{h}\cdot\sigma(\tau_{\tilde{h}^{-1}}(S))\cdot h.\]
Now we finally can prove Theorem~\ref{thm-cogredpar} by standard fibre bundle theory. If $L=\{L_\lambda\}_{\lambda\in I}$ is a~path in $\mathcal{FS}^i(H)^G$, then $L(I)$ is contained in a path component of $\mathcal{FS}^i(H)^G$, which we henceforth call~$C$. Let $S\in C$ be arbitrary, $Q$ the associated symmetry as in the first step and note that $C\subset B_{Q}$ by the second step. Now we consider the pullback $(E,I,\pi)$ of the fibre bundle
\begin{align}\label{piQ_S}
\pi_{Q} \colon \ \GL(H)^G\times\mathcal{KS}(H)^G\rightarrow B_{Q}
\end{align}
by $L$. This bundle has the set
\[E=\bigl\{(\lambda,h)\in I\times \bigl(\GL(H)^G\times\mathcal{KS}(H)^G\bigr) \mid L_\lambda=\pi_Q(h)\bigr\}\]
as total space and the bundle projection $\pi$ is the restriction of the projection onto the first component. As $I$ is contractible, $(E,I,\pi)$ is trivial and consequently has a globally defined section. Now the projection onto the second component $I\times\bigl(\GL(H)^G\times\mathcal{KS}(H)^G\bigr)\rightarrow\GL(H)^G\times\mathcal{KS}(H)^G$ is a bundle map from $E$ to the total space of \eqref{piQ_S}. By composing with this map, sections of $(E,I,\pi)$ yield liftings of $L$. Thus, we~have found the desired map $\tilde{L} \colon I\rightarrow\GL(H)^G\times\mathcal{KS}(H)^G$ such that $L_\lambda=\pi_Q\circ\tilde{L}_\lambda$ for all $\lambda\in I$ and so Theorem~\ref{thm-cogredpar} is proved.
\end{proof}

\subsection{\texorpdfstring{$\boldsymbol{G}$}{G}-equivariant approximation schemes}
In the remainder of the proof, we need finite-dimensional approximations as constructed for the classical spectral flow in \cite{JacoboUniqueness} by an orthonormal basis of the Hilbert space~$H$. The following lemma adapts the construction to the $G$-equivariant setting.

\begin{Lemma}\label{lem-Zorn}
For every symmetry $Q$, there is an increasing sequence of finite-dimensional $G$\nobreakdash-in\-vari\-ant subspaces $H_n\subset H$, $n\in\mathbb{N}$, such that the orthogonal projections $P_n$ onto $H_n$ commute with $Q$ and weakly converge to the identity.
\end{Lemma}

\begin{proof}
We first consider the case that $Q=I_H$ or $Q=-I_H$ and we aim to use Zorn's lemma. Let $\mathcal{S}$ be the set of all subsets $S$ of $H$ such that
\begin{itemize}\itemsep=0pt
 \item[(i)] $\langle x,y\rangle=\delta_{xy}$ for all $x,y\in S$,
 \item[(ii)] for all $x\in S$ there is a subspace $V\subset H$ of finite dimension that contains the orbit of $x$, i.e.,
$Gx:=\{gx \mid g \in G\}\subset V$.
\end{itemize}
We first note that $\mathcal{S}\neq\varnothing$ as every representation of $G$ on an infinite-dimensional Banach space has a finite-dimensional subrepresentation by \cite[Corollary~5.4\,(a)]{IzeVignoli}, which uses Zorn's lemma as well. Now $\mathcal{S}$ can be partially ordered by inclusion, and for any totally ordered $\mathcal{T}\subset\mathcal{S}$, the union of all elements in $\mathcal{T}$ satisfies (i) and (ii) and thus is an upper bound for $\mathcal{T}$. Consequently, by Zorn's lemma there exists a maximal element $S^\ast$ of $\mathcal{S}$, which moreover is countable as $S^\ast$ is orthonormal and $H$ is separable. We set $S^\ast=\{e_1,e_2,\dots\}$ and let $H_n$ be the intersection of all $G$-invariant subspaces of $H$ that contain $\{e_1,\dots,e_n\}$, which is of finite dimension by (ii). Moreover, the union $U:=\overline{\bigcup^\infty_{n=1}H_n}$ is a $G$-invariant subspace of $H$.

We now assume by contradiction that $U\neq H$. Then $U^\perp$ is a $G$-invariant subspace and thus contains a finite-dimensional subrepresentation, which again follows by \cite[Corollary~5.4\,(a)]{IzeVignoli}. Let $v$ be an element of norm $1$ of this finite-dimensional subrepresentation of $U^\perp$. Then \mbox{$S^\ast\cup\{v\}\in\mathcal{S}$} is larger than $S^\ast$ in contradiction to the maximality. Thus, $U=H$, and this in particular implies that $(P_n)_{n\in\mathbb{N}}$ weakly converges to the identity.

We now consider a symmetry $Q=P-(I_H-P)$ for some orthogonal projection $P$ that is neither $0$ nor the identity. Let $H^+$ be the image of $P$ and $H^-$ its kernel. As~$Q$ is~$I_{H^+}$ on~$H^+$ and $-I_{H^-}$ on $H^-$, by the argument from above, there are finite-dimensional $G$-invariant subspaces $H^\pm_n\subset H^\pm_{n+1}\subset H^\pm$, $n\in\mathbb{N}$, such that $P^\pm_nu\xrightarrow{n\rightarrow\infty} u$, $u\in H^\pm$, where $P^\pm_n$ are the orthogonal projections onto $H^\pm_n$ in $H^\pm$. If we set $H_n:=H^+_n\oplus H^-_n$, then $P_n:=P^+_n+P^-_n$ is the orthogonal projection onto $H_n$ if we regard $P^\pm_n$ as orthogonal projection on $H$ extended to $H^\mp$ by zero. Now~$H_n$~is an increasing sequence of invariant spaces and $(P_n)_{n\in\mathbb{N}}$ weakly converges to the identity. The projections $P_n$ also commute with $Q$ as
\begin{align*}
P_nQ-QP_n&=(P^+_n+P^-_n)(P-(I_H-P))-(P-(I_H-P))(P^+_n+P^-_n)\\
&=P^+_nP-P^-_n(I_H-P)-PP^+_n+(I_H-P)P^-_n=0,
\end{align*}
where we have used that $P_1P_2=P_2P_1=P_1$ if $\im(P_1)\subset\im(P_2)$ for any orthogonal projections $P_1$, $P_2$.
\end{proof}

Note that as $P_n$ commutes with $Q$ by the previous lemma, it follows that $Q(H_n)=H_n$ as well as $Q(H^\perp_n)=H^\perp_n$.

\subsection{Invariance under conjugation}\label{sect-OG}
We show the following further property of $\mu$, which is of independent interest. In what follows, we denote by $\mathcal{O}(H)$ the orthogonal operators on the real Hilbert space $H$.

\begin{Proposition}\label{prop-OG}
If $\mu \colon \Omega\bigl(\mathcal{FS}(H)^G,\GL(H)\cap\mathcal{FS}(H)^G\bigr)\rightarrow \mathrm{RO}(G)$ has the properties $(\mathcal{Z}_G)$, $(\mathcal{A}_G)$, $(\mathcal{H}_G)$ and $(\mathcal{M}_G)$, then
\begin{align}\label{OG}
\mu(U^\ast LU)=\mu(L)\in \mathrm{RO}(G)
\end{align}
for any $U\in\mathcal{O}(H)^G$ and $L\in\Omega\bigl(\mathcal{FS}(H)^G,\GL(H)\cap\mathcal{FS}(H)^G\bigr)$, where $U^\ast LU$ stands for the path $\{U^\ast L_\lambda U\}_{\lambda\in I}$.
\end{Proposition}

\begin{proof}
Let $R$ denote the at most countable set of all isomorphism classes of irreducible representations of $G$. If $H_\nu$ denotes the isotypical component of $H$ modelled on $\nu\in R$, then we have the decomposition
\[H=\bigoplus_{\nu\in R_H}H_\nu,\]
where $R_H\subset R$ is the subset of representations $\nu$ such that $H_\nu\neq\{0\}$. We now claim that for any path $L=\{L_\lambda\}_{\lambda\in I}$ in $\mathcal{FS}(H)^G$, there is a finite subset $R'\subset R_H$ such that
\begin{align}\label{lem-splitting}
H=\Biggl(\bigoplus_{\nu\in R'}H_\nu\Biggr)\oplus\Biggl(\bigoplus_{\nu\in R_H\setminus R'}H_\nu\Biggr)=:V\oplus W,
\end{align}
$L_\lambda(V)\subset V$, $L_\lambda(W)\subset W$, and $L_\lambda|_W \colon W\rightarrow W$ is an isomorphism for all $\lambda\in I$.

To show \eqref{lem-splitting}, let $0=\lambda_0<\lambda_1<\dots<\lambda_N=1$ and $a_i>0$ be as in the definition \eqref{sfl-equiv} of the equivariant spectral flow, which in particular means that
\[[\lambda_{i-1},\lambda_i]\ni\lambda\mapsto\chi_{[-a_i,a_i]}(L_\lambda)\in\mathcal{L}(H)^G\]
is a continuous family of finite rank projections for $i=1,\dots,N$. It follows from \cite[Lemma~2.2]{MJN21} that any $\im(\chi_{[-a_i,a_i]}(L_\lambda))$, $\lambda\in[\lambda_{i-1},\lambda_i]$, and the space $\im(\chi_{[-a_i,a_i]}(L_{\lambda_i}))$ for $\lambda=\lambda_i$ are isomorphic $G$-representations. If we now set $E^i:=\im(\chi_{[-a_i,a_i]}(L_{\lambda_i}))$, then this is a $G$\nobreakdash-rep\-re\-sen\-ta\-tion of finite dimension and thus
\[E^i=\bigoplus_{\nu\in R_i}E^i_\nu\]
for some finite subset $R_i\subset R_H$, where again $E^i_\nu$ is the isotypical component of $E^i$ modelled on $\nu$. For $R':=R_1\cup\dots\cup R_N$, the set $R'$ is finite and any $\im(\chi_{[-a_i,a_i]}(L_{\lambda}))$, $\lambda\in[\lambda_{i-1},\lambda_i]$, $i=1,\dots,N$, is a subrepresentation of
\[V:=\bigoplus_{\nu\in R'}H_\nu.\]
Consequently, $L_\lambda|_W \colon W\rightarrow W$, $\lambda\in I$, is an isomorphism for
\[W:=\bigoplus_{\nu\in R_H\setminus R'}H_\nu\]
and \eqref{lem-splitting} is shown.

Next, let us recall that any isotypical component $H_\nu$ is $G$-isomorphic to $E\otimes \nu$ for some Hilbert space $E$ on which $G$ acts trivially, and $\mathcal{L}(H_\nu)^G$ is isomorphic to $\mathcal{L}(E)\otimes\mathbb{K}$ (cf.\ \cite{Broecker}), where $\mathbb{K}=\mathbb{R}$, $\mathbb{K}=\mathbb{C}$, or $\mathbb{K}=\mathbb{H}$. Thus, it follows from Kuiper's theorem that $\mathcal{O}(H_\nu)^G$ is contractible if the isotypical component $H_\nu$ is of infinite dimension. Thus, $\mathcal{O}(H_\nu)^G$ is in particular path-connected in this case.

To finally show Proposition~\ref{prop-OG}, we let $H=V\oplus W$ as in \eqref{lem-splitting} and $U\in\mathcal{O}(H)^G$. Then
\begin{align*}
\mu(U^\ast LU)=\mu(U^\ast L|_V U)+\mu(U^\ast L|_WU)=\mu(U^\ast L|_V U),
\end{align*}
where we have used $(\mathcal{A}_G)$, $(\mathcal{Z}_G)$ and that $L_\lambda|_W$ is invertible for all $\lambda\in I$. We set
\[R_1:=\{\mu\in R' \mid \dim(H_\nu)=\infty\},\qquad R_2:=R'\setminus R_1,\]
where
\[V:=\bigoplus_{\nu\in R'}H_\nu\]
for some finite set $R'$ of representations of $G$ as in \eqref{lem-splitting}. Now
\begin{align*}
V=\Biggl(\bigoplus_{\nu\in R_1}H_\nu\Biggr)\oplus\Biggl(\bigoplus_{\nu\in R_2}H_\nu\Biggr),
\end{align*}
where both sums are finite and we obtain from $(\mathcal{A}_G)$
\[\mu(U^\ast L|_VU)=\sum_{\nu\in R_1}\mu(U^\ast L|_{H_\nu}U)+\sum_{\nu\in R_2}\mu(U^\ast L|_{H_\nu}U)\in \mathrm{RO}(G).\]
We have for $\nu\in R_1$
\[\mu(U^\ast L|_{H_\nu}U)=\mu(L|_{H_\nu})\]
by $(\mathcal{H}_G)$ as $\mathcal{O}(H_\nu)^G$ is path-connected. Moreover, we obtain by $(\mathcal{M}_G)$ for $\nu\in R_2$
\begin{align*}
\mu(U^\ast L|_{H_\nu}U)&=[E^-(U^\ast L_0|_{H_\nu}U)]-[E^-(U^\ast L_1|_{H_\nu}U)]\\
&=[E^-(L_0|_{H_\nu})]-[E^-( L_1|_{H_\nu})]=\mu(L|_{H_\nu}),
\end{align*}
where we use that $U^\ast$ is a $G$-equivariant isomorphism between $E^-(L_\lambda|_{H_\nu})$ and $E^-(U^\ast L_\lambda|_{H_\nu}U)$ for $\lambda\in I$. In summary,
\begin{align*}
\mu(U^\ast LU)&=\mu(U^\ast L|_V U)=\sum_{\nu\in R_1}\mu(U^\ast L|_{H_\nu}U)+\sum_{\nu\in R_2}\mu(U^\ast L|_{H_\nu}U)\\
&=\sum_{\nu\in R_1}\mu(L|_{H_\nu})+\sum_{\nu\in R_2}\mu(L|_{H_\nu})=\mu(L|_V)=\mu(L|_V)+\mu(L|_W)=\mu(L),
\end{align*}
where we have used $(\mathcal{A}_G)$ and $(\mathcal{Z}_G)$ in the equalities of the second line. Consequently, Proposition~\ref{prop-OG} is shown.
\end{proof}

\begin{Remark}
We noted below Corollary~\ref{cor-uniqueness} that Theorem~\ref{thm-uniqueness} also holds if we consider the $G$\nobreakdash-equiv\-ari\-ant spectral flow for complex Hilbert spaces, and the proof simplifies in this case. The reason for the latter is that the set of $G$-equivariant unitary operators $\mathcal{U}(H)^G$ is path-connected (cf.\ \cite[Section~2.2]{shortdegree} and \cite{Broecker}), so that a complex version of~\eqref{OG} follows from $(\mathcal{H}_G)$ by connecting~$U$ to the identity.
\end{Remark}

\subsection{Finite-dimensional approximation and end of proof}

We now consider a path $L\in\Omega\bigl(\mathcal{FS}(H)^G,\GL(H)\cap\mathcal{FS}(H)^G\bigr)$. By Theorems~\ref{thm-cogrpm} and~\ref{thm-cogredpar}, there is a $G$-equivariant cogredient parametrix for $L$, i.e., a path $M=\{M_\lambda\}_{\lambda\in I}$ in $\GL(H)^G$ such that
\begin{align*}
M^\ast_\lambda L_\lambda M_\lambda= Q+K_\lambda,\qquad \lambda\in I,
\end{align*}
where $K_\lambda$ is $G$-equivariant, selfadjoint and compact, and $Q$ is a $G$-equivariant symmetry. We set $\tilde{L}_\lambda=M^\ast_\lambda L_\lambda M_\lambda$ and consider the path $\tilde{L}=\{\tilde{L}_\lambda\}_{\lambda\in I}$, which is an element of $\Omega\bigl(\mathcal{FS}(H)^G,\GL(H)\cap\mathcal{FS}(H)^G\bigr)$ so that $\mu\bigl(\tilde{L}\bigr)$ is defined.

\begin{Lemma}
The paths of operators $L$ and $\tilde{L}$ from above have the same $\mu$-class in $\mathrm{RO}(G)$, i.e.,
\[\mu(L)=\mu\bigl(\tilde{L}\bigr)\in \mathrm{RO}(G).\]
\end{Lemma}

\begin{proof}
The path $\tilde{L}$ is homotopic to $\{M^\ast_0L_\lambda M_0\}_{\lambda\in I}$ and the corresponding homotopy does not affect $\mu$ by $(\mathcal{H}_G)$ as $L_0$ and $L_1$ are invertible. Next, we consider the polar decomposition $M_0=UR$ of $M_0$ into two $G$-equivariant operators which are given by $U=M_0(M^\ast_0M_0)^{-\frac{1}{2}}$ and $R=(M^\ast_0M_0)^\frac{1}{2}$. Note that $U$ is orthogonal and $R$ is selfadjoint and positive. Then we have
\[M^\ast_0L_\lambda M_0=RU^\ast L_\lambda UR,\qquad\lambda\in I,\]
and see that the homotopy
\[\{((1-s)R+sI_H)U^\ast L_\lambda U((1-s)R+sI_H)\}_{(s,\lambda)\in I\times I}\]
deforms $\{M^\ast_0L_\lambda M_0\}_{\lambda\in I}$ into the path $\{U^\ast L_\lambda U\}_{\lambda\in I}$. As $(1-s)R+sI_H\in\GL(H)^G$ for all $s\in I$, again this homotopy does not affect $\mu$ by $(\mathcal{H}_G)$. Finally, the claim follows from Pro\-po\-si\-tion~\ref{prop-OG}.
\end{proof}

For simplicity of notation, we henceforth assume that $L_\lambda=Q+K_\lambda$, $\lambda\in I$, where the operators $K_\lambda$ are compact, selfadjoint and $G$-equivariant, and $Q=P-(I_H-P)$ for some $G$\nobreakdash-equiv\-ari\-ant orthogonal projection $P$. By Lemma~\ref{lem-Zorn}, there is an approximation scheme for $Q$, i.e., an~increasing sequence $H_n$, $n\in\mathbb{N}$, of finite-dimensional $G$-invariant subspaces of $H$ such that the orthogonal projections $P_n$ onto $H_n$ commute with $Q$ and weakly converge to the identity. We~now consider for $n\in\mathbb{N}$ the paths $P_nLP_n=\{P_nL_\lambda P_n\}_{\lambda\in I}$ made by the $G$-equivariant operators $P_nL_\lambda P_n \colon H_n\rightarrow H_n$.

\begin{Lemma}
There is $n_0\in\mathbb{N}$ such that $P_n LP_n$ has invertible endpoints for all $n\geq n_0$ and
\[\mu(L)=\mu(P_n LP_n)\in \mathrm{RO}(G).\]
\end{Lemma}

\begin{proof}
We first consider the operator
\[
(I_H-P_n)L_\lambda|_{H^\perp_n}=Q+(I_H-P_n)K_\lambda|_{H^\perp_n}
\]
 and aim to show that it is invertible for sufficiently large $n\in\mathbb{N}$. As it is a compact perturbation of an invertible operator, it is a Fredholm operator of index $0$ and so we only need to show the injectivity.

Since the family $\{K_\lambda\}_{\lambda\in I}$ consists of compact operators, the subset $\{K_\lambda(u) \mid \lambda \in I,\, \|u\|=1\}$ of $H$ is relatively compact. As $P_n$ weakly converges to the identity, $I_H-P_n$ uniformly converges to $0$ on compact subsets of $H$. Consequently, there exists $n_0\in\mathbb{N}$ such that
\[\|(I_H-P_n)K_\lambda u\|\leq\frac{1}{2}\|u\|,\qquad u\in H, \quad \lambda\in I,\quad n\geq n_0.\]
As $Q$ is orthogonal, we have $\|Qu\|=\|u\|$, $u\in H$, and thus
\begin{align*}
\|(I_H-P_n)L_\lambda u\|=\|Qu+(I_H-P_n)K_\lambda u\|\geq\frac{1}{2}\|u\|,\qquad u\in H^\perp_n,
\end{align*}
which indeed shows that $(I_H-P_n)L_\lambda|_{H^\perp_n}$ is injective for $n\geq n_0$.

As second step, we consider the homotopy
\[
h(s,\lambda)=sL_\lambda+(1-s)((I_H-P_n)L_\lambda(I_H-P_n)+P_nL_\lambda P_n),\qquad (s,\lambda)\in I\times I,
\]
in $\mathcal{FS}(H)^G$ and aim to show that the endpoints are invertible for sufficiently large $n\in\mathbb{N}$, i.e., $h(s,0), h(s,1)\in \GL(H)$ for all $s\in I$, $n\geq n_0$, and some $n_0\in\mathbb{N}$. We first note that
\begin{align*}
h(s,\lambda)=Q+sK_\lambda+(1-s)((I_H-P_n)K_\lambda(I_H-P_n)+P_nK_\lambda P_n),
\end{align*}
and thus $h(s,\lambda)$ are indeed all Fredholm operators of index $0$. Let us now assume by contradiction that an $n_0$ as claimed does not exist. Then there are a sequence $(u_n)_{n\in\mathbb{N}}$, $\|u_n\|=1$, in $H$ and a~sequence $(s_n)_{n\in\mathbb{N}}$ in $I$ such that
\[Qu_n+s_nK_0u_n+(1-s_n)((I_H-P_n)K_0(I_H-P_n)u_n+P_nK_0 P_nu_n)=0,\qquad n\in\mathbb{N}.\]
Now $K_0$ is compact and $P_n$ converges on compact subsets of $H$ to the identity, which implies that there is a subsequence $(u_{n_j})_{j\in\mathbb{N}}$ such that $(Qu_{n_j})_{j\in\mathbb{N}}$ converges. Without loss of generality, we can assume that $(s_{n_j})$ converges to some $s^\ast\in I$. Then the invertibility of $Q$ shows that $(u_{n_j})_{j\in\mathbb{N}}$ converges to some $u\in H$ of norm $1$. Thus,
\[\lim_{n\rightarrow\infty} (I_H-P_n)K_0(I_H-P_n)u_n=0,\qquad \lim_{n\rightarrow\infty}P_nK_0 P_n u_n=K_0 u,\]
and so
\[L_0 u=Qu+K_0u=Qu+s^\ast K_0 u+(1-s^\ast)K_0 u=0,\]
which contradicts the invertibility of $L_0$. Of course, the same argument applies to the invertible operator $L_1$.

Thus, we have shown that there is some $n_0\in\mathbb{N}$ such that for $n\geq n_0$, $h(s,0)$ and $h(s,1)$ are invertible for all $s\in I$ and the operators $(I_H-P_n)L_\lambda|_{H^\perp_n}$, $\lambda\in I$, are invertible. This in particular shows that $P_nL_0 P_n$ and $P_n L_1 P_n$ are invertible. Finally, we obtain from $(\mathcal{H}_G)$ and $(\mathcal{A}_G)$ that for $n\geq n_0$
\[\mu(L)=\mu(h(1,\cdot))=\mu(h(0,\cdot))=\mu((I_H-P_n)L(I_H-P_n)+P_nL P_n)=\mu(P_n LP_n),\]
as claimed.
\end{proof}

Now the end of the proof is rather simple. We first note that by $(\mathcal{M}_G)$ and the previous lemma
\begin{align}\label{eqref-finreduction}
\mu(L)=\mu(P_n LP_n)=[E^-(P_n L_0 P_n)]-[E^-(P_n L_1 P_n)]\in \mathrm{RO}(G)
\end{align}
for sufficiently large $n\in\mathbb{N}$. As the right-hand side does not depend on $\mu$ but only on the path~$L$, any two choices of a map $\mu$ ultimately must be the same by the equalities in \eqref{eqref-finreduction}.

\section[The G-equivariant spectral flow as G-equivariant Maslov index]{The $\boldsymbol{G}$-equivariant spectral flow\\ as $\boldsymbol{G}$-equivariant Maslov index}

\subsection{The \texorpdfstring{$\boldsymbol{G}$}{G}-equivariant Maslov index}
To introduce the $G$-equivariant Maslov index, we firstly need a slightly more general definition of the $G$-equivariant spectral flow. As selfadjoint operators are closed, the set of (generally) unbounded selfadjoint Fredholm operators $\mathcal{A}\colon \mathcal{D}(\mathcal{A})\subset H\rightarrow H$ naturally is a metric space by
\begin{align}\label{gap}
d(\mathcal{A}_1,\mathcal{A}_2)=\|P_1-P_2\|,
\end{align}
where $P_i$, $i=1,2$, denote the orthogonal projections onto the graphs of $\mathcal{A}_i$ in $H\times H$. Thus, it is sensible to consider paths $\mathcal{A}_\lambda \colon \mathcal{D}(\mathcal{A}_\lambda)\subset H\rightarrow H$ of selfadjoint Fredholm operators, where now the domain $\mathcal{D}(\mathcal{A}_\lambda)$ depends on the parameter $\lambda\in I$. It was shown by Boo{\ss}-Bavnbek, Lesch and Phillips in \cite{UnbSpecFlow} that the spectral flow can be defined verbatim by the same formula \eqref{sfl} and that all its basic properties still hold in this more general setting. Now, if in addition $G$ is a compact Lie group acting orthogonally on $H$, $\mathcal{D}(\mathcal{A}_\lambda)$ is $G$-invariant and $\mathcal{A}_\lambda$ is $G$-equivariant for all $\lambda\in I$, it should come as no surprise that \eqref{sfl-equiv} defines a $G$-equivariant spectral flow $\sfl_G(\mathcal{A})$ in $\mathrm{RO}(G)$ which in particular satisfies corresponding properties $(\mathcal{Z}_G)$, $(\mathcal{A}_G)$ and $(\mathcal{H}_G)$ in this setting. Actually, the $G$-equivariant spectral flow has already been introduced for paths of unbounded operators in \cite{MJN21}, where constant domains were only assumed for the sake of simplicity.

Let us next recall that a symplectic Hilbert space is a real Hilbert space $E$ with an invertible bounded operator $J \colon E\rightarrow E$ such that $J^\ast=-J$ and $J^2=-I_E$. This yields a symplectic form on $E$ by $\omega(u,v)=\langle Ju,v\rangle$. A closed subspace $L\subset E$ is called Lagrangian if $L$ and its annihilator $\{u\in E \mid \omega(u,v)=0\,\,\text{for all}\,\, v\in L\}$ coincide, which is equivalent to
\begin{align}\label{Lag}
L^\perp=J(L).
\end{align}
In the latter formula, the orthogonal complement is meant with respect to the scalar product of the Hilbert space $E$ (cf. \cite{Furutani}). The set $\Lambda(E,\omega)$ of all Lagrangian subspaces of $E$ naturally is a~metric space by
\begin{align}\label{gapII}
d(L_1,L_2)=\|P_1-P_2\|,
\end{align}
where $P_i$, $i=1,2$, are the orthogonal projections onto the closed subspaces $L_i$. Let us recall that two closed subspaces $U$, $W$ are called a Fredholm pair if $\dim(U\cap W)$ and $\codim(U+W)$ are of finite dimension. We now fix some $W\in\Lambda(E,\omega)$ and set
\[\mathcal{FL}_W(E,\omega)=\{L\in \Lambda(E,\omega) \mid (L,W) \ \text{is a Fredholm pair}\}.\]
Boo{\ss}-Bavnbek and Furutani constructed in \cite{Bernhelm} a Maslov index $\mu_{\rm Mas}(\Gamma,W)\in\mathbb{Z}$ for paths $\Gamma=\{\Gamma_\lambda\}_{\lambda\in I}$ in $\mathcal{FL}_W(E,\omega)$. Roughly speaking, the Maslov index counts the intersections of $\Gamma_\lambda$ with~$W$ including dimensions whilst the parameter $\lambda$ is travelling along the interval $I$. If~$E$~is of finite dimension, then every pair of subspaces is Fredholm and this construction is the classical one \cite{Robbin-SalamonMAS}. A rigorous construction of the Maslov index in symplectic Hilbert spaces is rather involved (see \cite{Furutani}) and uses a spectral flow for paths of unitary operators. In the finite-dimensional case, Capell, Lee and Miller showed in a seminal paper~\cite{Cappell} that it can equivalently be defined as the spectral flow of a path of selfadjoint differential operators (see also~\cite{IJW}), which can make it very accessible. This link between the Maslov index and the spectral flow was generalised to symplectic Hilbert spaces by the last author in~\cite{ODE}, and now paves the way for introducing a $G$-equivariant Maslov index. Let us note that the more general case of relative Maslov indices was considered in~\cite{Furutani} as well and the below construction can be generalised along the same lines.

Let $E$ be a symplectic Hilbert space with symplectic form $\omega$ represented by the bounded invertible operator $J \colon E\rightarrow E$. We assume that $G$ is a compact Lie group acting orthogonally and symplectically on $E$, i.e., $g^\ast g=g g^\ast=I_E$ and $g^\ast Jg=J$ for all $g\in G$. Let $W\in\Lambda(E,\omega)$ be $G$-invariant and $\Gamma=\{\Gamma_\lambda\}_{\lambda\in I}$ a path of $G$-invariant spaces in $\mathcal{FL}_W(E,\omega)$. We define differential operators on $L^2(I,E)$ by
\begin{align}\label{Aop}
\mathcal{A}^\Gamma_\lambda: \ \mathcal{D}\bigl(\mathcal{A}^\Gamma_\lambda\bigr)\subset L^2(I,E)\rightarrow L^2(I,E),\qquad \mathcal{A}^\Gamma_\lambda u=Ju',
\end{align}
where
\[\mathcal{D}(\mathcal{A}^\Gamma_\lambda)=\bigl\{u\in H^1(I,E) \mid u(0) \in\Gamma_\lambda,\, u(1)\in W\bigr\}.\]
It was shown in \cite{ODE} that $\mathcal{A}^\Gamma=\bigl\{\mathcal{A}^\Gamma_\lambda\bigr\}_{\lambda\in I}$ is a continuous path of selfadjoint Fredholm operators with respect to \eqref{gap} if $\Gamma=\{\Gamma_\lambda\}_{\lambda\in I}$ is continuous with respect to \eqref{gapII}. Note that the dimension of the kernel of $\mathcal{A}^\Gamma_\lambda$ exactly is the dimension of the intersection $\Gamma_\lambda\cap W$. Theorem~B of \cite{ODE} states that
\[\mu_{\rm Mas}(\Gamma,W)=\sfl\bigl(\mathcal{A}^\Gamma\bigr)\in\mathbb{Z}.\]
Note that the orthogonal action of $G$ on $E$ naturally transfers to an orthogonal action on $L^2(I,E)$, $\mathcal{D}\bigl(\mathcal{A}^\Gamma_\lambda\bigr)$ is invariant under this action and, since the action is symplectic, $\mathcal{A}^\Gamma_\lambda$ is $G$\nobreakdash-equiv\-ari\-ant. Thus, the $G$-equivariant spectral flow of $\mathcal{A}^\Gamma$ is defined as explained above, and we can now make the following definition.

\begin{Definition}\label{def-Maslov}
Let $E$ be a symplectic Hilbert space with symplectic form $\omega$ represented by $J \colon E\rightarrow E$ and $G$ a compact Lie group acting orthogonally and symplectically on $E$. If~$W\in\Lambda(E,\omega)$ is $G$-invariant and $\Gamma=\{\Gamma_\lambda\}_{\lambda\in I}$ is a path in $\mathcal{FL}_W(E,\omega)$ such that each $\Gamma_\lambda$ is $G$-invariant, then the $G$-equivariant Maslov index is defined by
\[\mu^G_{\rm Mas}(\Gamma,W)=\sfl_G\bigl(\mathcal{A}^\Gamma\bigr)\in \mathrm{RO}(G),\]
where $\mathcal{A}^\Gamma$ is the associated path of selfadjoint Fredholm operators in \eqref{Aop}.
\end{Definition}

\subsection{The \texorpdfstring{$\boldsymbol{G}$}{G}-equivariant spectral flow as \texorpdfstring{$\boldsymbol{G}$}{G}-equivariant Maslov index}
Our aim is now to show an example in which we compute the $G$-equivariant Maslov index by Theorem~\ref{thm-uniqueness}. We let once again $H$ be a real separable Hilbert space and $G$ a compact Lie group which acts orthogonally on $H$. Moreover, we assume that $L=\{L_\lambda\}_{\lambda\in I}$ is a path in $\mathcal{FS}(H)^G$ as in Theorem~\ref{thm-uniqueness}. Note that the space $E:=H\times H$ is a symplectic Hilbert space with respect to the standard symplectic form, which is represented by $J \colon E\rightarrow E$,
\[J=\begin{pmatrix}
0&-I_H\\
I_H&0
\end{pmatrix}.\]
We now consider the subspaces $W=H\times\{0\}\subset E$ and
\[\Gamma_\lambda=\gra(L_\lambda)=\{(u,L_\lambda u)\mid u\in H\},\qquad \lambda\in I.\]
The following lemma shows that the classical Maslov index $\mu_{\rm Mas}(\Gamma,W)$ of $\Gamma$ with respect to $W$ is defined.

\begin{Lemma}\label{lemma-FLW}
For $W$ as above, $\Gamma=\{\Gamma_\lambda\}_{\lambda\in I}$ is a path in $\mathcal{FL}_W(E,\omega)$.
\end{Lemma}

\begin{proof}
We first note that it follows from \eqref{Lag} that $W$ and $\Gamma_\lambda$, $\lambda\in I$, are Lagrangian. Moreover, the Fredholm property of $L_\lambda$ implies that $\Gamma_\lambda\in\mathcal{FL}_W(E,\omega)\subset\Lambda(E,\omega)$. Finally, every path $L=\{L_\lambda\}_{\lambda\in I}$ in $\mathcal{FS}(H)$ is continuous with respect to the metric \eqref{gap} (cf.\ \cite{Lesch}) and thus the path $\Gamma$ is continuous in $\Lambda(E,\omega)$ with respect to \eqref{gapII} by definition.
\end{proof}

We now consider the diagonal action $g(u,v)=(gu,gv)$ on $E=H\times H$ and note that this action is orthogonal and symplectic. Moreover, $W$ and $\gra(L_\lambda)$ are $G$-invariant, where we use in the latter case that $L_\lambda$ is $G$-equivariant. Thus, the $G$-equivariant Maslov index is defined.

\begin{Theorem}\label{thm-Maslov}
Let $G$ be a compact Lie group that acts orthogonally on $H$. For every path $L=\{L_\lambda\}_{\lambda\in I}$ in $\mathcal{FS}(H)^G$ with invertible ends
\[\sfl_G(L)=\mu^G_{\rm Mas}(\Gamma,W)\in \mathrm{RO}(G),\]
where $\Gamma=\{\gra(L_\lambda)\}_{\lambda\in I}$ and $W=H\times\{0\}$.
\end{Theorem}

\begin{proof}
We aim to use Theorem~\ref{thm-uniqueness}, thus set
\[\mu \colon \ \Omega\bigl(\mathcal{FS}(H)^G,\GL(H)\cap\mathcal{FS}(H)^G\bigr)\rightarrow \mathrm{RO}(G),\qquad \mu(L)=\mu^G_{\rm Mas}(\Gamma,W),\]
and now show $(\mathcal{Z}_G)$, $(\mathcal{A}_G)$, $(\mathcal{H}_G)$ and $(\mathcal{M}_G)$ for $\mu$. Let us emphasize that we can use these properties for the spectral flow of the path $\mathcal{A}^\Gamma$ in Definition~\ref{def-Maslov}, which we henceforth do without further reference.

We first consider $(\mathcal{Z}_G)$ for a path $L=\{L_\lambda\}_{\lambda\in I}$. If $L_\lambda$ is invertible for all $\lambda\in I$, then\linebreak $\ker(L_\lambda)\cap W=\{0\}$ and thus $\mathcal{A}^\Gamma_\lambda$ is invertible for all $\lambda\in I$ as it in particular is a Fredholm operator of index $0$. Consequently, $\mu(L)=0$ by Definition~\ref{def-Maslov} and $(\mathcal{Z}_G)$ for the path $\mathcal{A}^\Gamma$.

For $(\mathcal{A}_G)$, we let $H=H_1\oplus H_2$ for two $G$-invariant subspaces that reduce the operators $L_\lambda$. Now $E=H\times H=(H_1\oplus H_2)\times (H_1\oplus H_2)$, and if we set $L^i_\lambda:=L_\lambda|_{H_i}$, $i=1,2$, then
\[\gra(L_\lambda)=\bigl\{\bigl((u,v),\bigl(L^1_\lambda u,L^2_\lambda v\bigr)\bigr)\in E \mid (u,v)\in H=H_1\oplus H_2\bigr\}.\]
We let $\mathcal{A}^\Gamma_\lambda$ be the operator for $\Gamma_\lambda=\gra(L_\lambda)$ in Definition~\ref{def-Maslov}. Moreover, for $E^i=H_i\times H_i$, $i=1,2$, we have
\[\gra(L^i_\lambda)=\bigl\{\bigl(u,L^i_\lambda u\bigr)\in E^i \mid u\in H^i\bigr\},\]
and denote by $\mathcal{A}^{\Gamma,i}_\lambda$ the corresponding operators in Definition~\ref{def-Maslov}. Let us now consider the operator $U \colon E\rightarrow E$ defined by $U(v_1,v_2,v_3,v_4)=(v_1,v_3,v_2,v_4)$. Then $U$ is orthogonal and $U\gra(L_\lambda)=\gra\bigl(L^1_\lambda\bigr)\oplus\gra\bigl(L^2_\lambda\bigr)$. Moreover, $U(\{(v_1,v_2,0,0)\in E \mid v_1\in H_1, v_2\in H_2\})=W_1\oplus W_2$, where $W_i=H_i\times\{0\}\subset E^i$, $i=1,2$. Now $U \colon
E\rightarrow E$ yields an orthogonal operator on $L^2(I,E)$, which we also denote by $U$, and we see that $U\mathcal{D}\bigl(\mathcal{A}^\Gamma_\lambda\bigr)=\mathcal{D}\bigl(\mathcal{A}^{\Gamma,1}_\lambda\oplus\mathcal{A}^{\Gamma,2}_\lambda\bigr)$. Moreover, a direct computation shows that
\[U^\ast JU=\begin{pmatrix}
0&-I_{H_1}&0&0\\
I_{H_1}&0&0&0\\
0&0&0&-I_{H_2}\\
0&0&I_{H_2}&0
\end{pmatrix}=\begin{pmatrix}
J_{E^1}&0\\
0&J_{E^2}
\end{pmatrix},\]
where $J_{E^i} \colon E^i\rightarrow E^i$, $i=1,2$, denotes the bounded invertible operator that represents the symplectic form in $E^i$. Consequently, $U^\ast\mathcal{A}^\Gamma_\lambda U=\mathcal{A}^{\Gamma,1}\oplus\mathcal{A}^{\Gamma,2}$, which shows that
\[\sfl_G\bigl(\mathcal{A}^\Gamma\bigr)=\sfl_G\bigl(\mathcal{A}^{\Gamma,1}\oplus\mathcal{A}^{\Gamma,2}\bigr) =\sfl_G\bigl(\mathcal{A}^{\Gamma,1}\bigr)+\sfl_G\bigl(\mathcal{A}^{\Gamma,2}\bigr)\in \mathrm{RO}(G)\]
by Proposition~\ref{prop-OG} and $(\mathcal{A}_G)$ for $\mathcal{A}^\Gamma$. Therefore, it follows from Definition~\ref{def-Maslov} that $\mu(L)=\mu(L|_{H_1}) +\mu(L|_{H_2})$, which is $(\mathcal{A}_G)$ for $L$.

To show $(\mathcal{H}_G)$ for $L$, we let $h\colon I\times I\rightarrow\mathcal{FS}(H)^G$ be a family such that $h(s,0)$ and $h(s,1)$ are invertible for all $s\in I$. As in the proof of Lemma~\ref{lemma-FLW}, the family $\Gamma=\{\gra(h(s,\lambda))\}_{(s,\lambda)\in I\times I}$ is continuous with respect to \eqref{gapII} and thus yields a homotopy $\mathcal{A}^\Gamma=\bigl\{\mathcal{A}^\Gamma_{(\lambda,s)}\bigr\}_{(\lambda,s)\in I\times I}$ of unbounded operators as in \eqref{Aop}. Since $(\mathcal{H}_G)$ holds for this homotopy, we obtain $\mu(h(0,\cdot))=\mu(h(1,\cdot))\in \mathrm{RO}(G)$ and thus $(\mathcal{H}_G)$ for $L$.

We now focus on $(\mathcal{M}_G)$ and assume that $H$ is of finite dimension. Let us consider $L_\lambda$ for a~fixed $\lambda\in I$, as well as the operator $\mathcal{A}^\Gamma_\lambda$ on $E=H\times H$. Note that $L_\lambda$ has a finite spectrum as~$H$~is of finite dimension, and $\mathcal{A}^\Gamma_\lambda$ has a discrete spectrum which only consists of eigenvalues of finite multiplicity as $\mathcal{A}^\Gamma_\lambda$ now has a compact resolvent (cf.\ \cite{Cappell}). The equation $Ju'=\alpha u$ has the fundamental solution
\[\cos(\alpha t)I_E-\sin(\alpha t)J,\qquad t\in I,\]
and thus we obtain all solutions that satisfy $u(0)=(u,L_\lambda u)\in\gra(L_\lambda)$ for some $u\in H$ by
\[u(t)=\begin{pmatrix}
\cos(\alpha t)u+\sin(\alpha t)L_\lambda u\\
-\sin(\alpha t)u+\cos(\alpha t) L_\lambda u
\end{pmatrix}.\]
As we also need to require that $u(1)\in H\times\{0\}$, we see that $u\in\mathcal{D}(\mathcal{A}^\Gamma_\lambda)$, if and only if $-\sin(\alpha) u+\cos(\alpha) L_\lambda u=0$. If $\alpha\in(-\frac{\pi}{2},\frac{\pi}{2})$, then the latter equation is equivalent to $L_\lambda u=\tan(\alpha) u$. Thus, the whole spectrum of $L_\lambda$ can be bijectively mapped to the eigenvalues of $\mathcal{A}_\lambda$ in $(-\frac{\pi}{2},\frac{\pi}{2})$. Now let $0=\lambda_0<\lambda_1<\dots<\lambda_N=1$ and $a_i\in(-\frac{\pi}{2},\frac{\pi}{2})$, $i=1,\dots,N$, be numbers as in the definition of the spectral flow in \eqref{sfl} for $\mathcal{A}^\Gamma$. As the operators $\mathcal{A}^\Gamma_\lambda$ have discrete spectrum, we can choose the $a_i$ sufficiently large, such that the whole spectrum of $L_\lambda$ is contained in $[-\tan(a_i),\tan(a_i)]$ for $\lambda_{i-1}\leq\lambda\leq\lambda_i$. This yields for any $a\geq\max\{\tan(a_1),\dots,\tan(a_N)\}$
\begin{align*}
\mu(L)&=\sfl_G(\mathcal{A}^\Gamma)=\sum^N_{i=1}{\bigl(\bigl[E\bigl(\mathcal{A}^\Gamma_{\lambda_i},[0,a_i]\bigr)\bigr] -\bigl[E\bigl(\mathcal{A}^\Gamma_{\lambda_{i-1}},[0,a_i]\bigr)\bigr]\bigr)}\\
&=\sum^N_{i=1}{[E(L_{\lambda_i},[0,a])]-[E(L_{\lambda_{i-1}},[0,a])]}=[E(L_{1},[0,a])]-[E(L_{0},[0,a])]\\
&=[E(L_{0},[-a,0))]-[E(L_{1},[-a,0))]=[E^-(L_0)]-[E^-(L_1)]\in \mathrm{RO}(G),
\end{align*}
where we have used in the third equality that mapping an eigenfunction $u$ of $\mathcal{A}^\Gamma_\lambda$ to the projection onto $H\times\{0\}$ of $u(0)$ is a $G$-equivariant isomorphism of eigenspaces. Moreover, the second last equality follows from
\begin{align*}
0&=[E(L_{1},[-a,a])]-[E(L_{0},[-a,a])]\\
&=([E(L_{1},[0,a])]-[E(L_{0},[0,a])])+([E(L_{1},[-a,0))]-[E(L_{0},[-a,0))]),
\end{align*}
which ends the proof of the theorem.
\end{proof}

Let us point out, without going into details, that it is also possible to prove
\[
\sfl_G(L)=\mu^G_{\rm Mas}(\Gamma,W)=\sfl_G\bigl(\mathcal{A}^\Gamma\bigr)\in \mathrm{RO}(G)
\] directly from the definition \eqref{sfl-equiv}. This requires a similar argument as for $(\mathcal{M}_G)$ in infinite dimensions as well as a careful choice of the partition $0=\lambda_0\leq\dots\leq\lambda_N=1$ and the numbers $a_i>0$, $i=1,\dots,N$, for both paths $L$ and $\mathcal{A}^\Gamma$.

\subsection{An example}
In this final part, we discuss a simple example of a path in the Fredholm Lagrangian Grassmannian of a Hilbert space for which the classical Maslov index vanishes but the $G$-equivariant Maslov index is non-trivial. Let $H$ be a real separable Hilbert space and $M=\{M_\lambda\}_{\lambda\in I}$ be a path in $\mathcal{FS}(H)$ with invertible endpoints. We consider the family $\Gamma_\lambda$ of subspaces of $E=H\oplus H\oplus H\oplus H$ given by
\[\Gamma_\lambda=\{(u,v,M_\lambda u,-M_\lambda v) \mid u,v \in H\}\qquad\text{and}\qquad W=H\oplus H\oplus\{0\}\oplus\{0\}.\]
As the operators $M_\lambda$ are selfadjoint and Fredholm, it is readily seen that actually $\Gamma=\{\Gamma_\lambda\}_{\lambda\in I}$ is a path in $\mathcal{FL}_W(E,\omega)$, where $\omega$ is represented by
\[J=\begin{pmatrix}
0&-I_{H\oplus H}\\
I_{H\oplus H}&0
\end{pmatrix}.\]
We now consider the $G=\mathbb{Z}_2$-action on $E$ given by
\begin{align}\label{action}
g(u_1,u_2,u_3,u_4)=(u_1,-u_2,u_3,-u_4),
\end{align}
where $g$ denotes the non-trivial element of $\mathbb{Z}_2$, and note that all spaces $\Gamma_\lambda$, $\lambda\in I$, as well as $W$ are invariant under this action. Thus, the $\mathbb{Z}_2$-equivariant Maslov index $\mu^G_{\rm Mas}(\Gamma,W)$ is defined as an element of $\mathrm{RO}(\mathbb{Z}_2)$. Now recall that all real irreducible representations of $\mathbb{Z}_2$ are one-dimensional and thus every real $k$-dimensional representation is up to isomorphism a $k\times k$ diagonal matrix of the form $\diag(1,\dots,1,-1,\dots,-1)$. Consequently, there is an isomorphism $\phi \colon \mathrm{RO}(\mathbb{Z}_2)\rightarrow\mathbb{Z}\oplus\mathbb{Z}$ of abelian groups given by{\samepage
\begin{align}\label{phi}
 \phi([E]-[F])=\bigl(\dim(E)-\dim(F),\dim\bigl(E^G\bigr)-\dim\bigl(F^G\bigr)\bigr),
\end{align}
where $E^G\subset E$ and $F^G\subset F$ denote the spaces of fixed points under the group action.}

The $G$-equivariant Maslov index $\mu^G_{\rm Mas}(\Gamma,W)$ can now be computed by Theorem~\ref{thm-Maslov}. Indeed, $\Gamma_\lambda$ is the graph of the selfadjoint Fredholm operator $L_\lambda\colon H\oplus H\rightarrow H\oplus H$ given by
\[L_\lambda=\begin{pmatrix}
M_\lambda&0\\
0&-M_\lambda
\end{pmatrix}.\]
Moreover, $L_\lambda$ is equivariant under the $\mathbb{Z}_2$-action $g(u,v)=(u,-v)$ on $H\oplus H$, where $g$ again denotes the non-trivial element of $\mathbb{Z}_2$. Note that \eqref{action} is the corresponding diagonal action and thus we are in the setting of Theorem~\ref{thm-Maslov}. It is not difficult to see from \eqref{phi} (see \cite[Section~2.4]{MJN21} for the full argument) that
\[\phi(\sfl_{\mathbb{Z}_2}(L))=(\sfl(L),\sfl(M))=(0,\sfl(M))\in\mathbb{Z}\oplus\mathbb{Z},\]
where we have used that the classical spectral flow \eqref{sfl} of $L$ has to vanish as the spectrum of each $L_\lambda$ is symmetric about $0$.

Thus, $\mu^{\mathbb{Z}_2}_{\rm Mas}(\Gamma,W)\in \mathrm{RO}(\mathbb{Z}_2)$ is non-trivial if and only if $\sfl(M)\in\mathbb{Z}$ is non-trivial. Finally, note that in case of the trivial group $G$, we see from Theorem~\ref{thm-Maslov} that
\[\mu_{\rm Mas}(\Gamma,W)=\sfl(L)=0\in\mathbb{Z},\]
i.e., the classical Maslov index~\cite{Bernhelm} vanishes.

\subsection*{Acknowledgements}
The authors were supported by the Deutsche Forschungsgemeinschaft (DFG, German Research Foundation) -- 459826435. Moreover, the authors wish to express their gratitude to the anonymous referees for their interest in our work and their valuable hints to improve the presentation of the paper.

\pdfbookmark[1]{References}{ref}
\LastPageEnding

\end{document}